\documentclass[a4paper,11pt]{article}

\usepackage[margin=0.6in]{geometry} % full-width

\usepackage{graphicx}% Include figure files

\usepackage{amsmath}
\usepackage{amsthm}
\usepackage{amssymb}

\usepackage[utf8]{inputenc}
\usepackage{hyperref}

\usepackage[sort&compress,numbers]{natbib}
\usepackage{subfigure}  % use for side-by-side figures
\usepackage{dcolumn}% Align table columns on decimal point
\usepackage{bm}% bold math
\usepackage{hyperref}
\hypersetup{colorlinks,allcolors=black}
\usepackage[mathlines]{lineno}% Enable numbering of text and display math
\newcommand{\beq}{\begin{equation}}
\newcommand{\eeq}{\end{equation}}

\newcommand{\ba}{\begin{array}}
\newcommand{\ea}{\end{array}}

\newcommand{\ben}{\begin{enumerate}}
\newcommand{\een}{\end{enumerate}}

\newtheorem{theorem}{Theorem}
\newcommand{\bth}{\begin{theorem}}

\newtheorem{definition}{definition}
\newcommand{\bdf}{\begin{definition}}
\newcommand{\edf}{\end{definition}}

\newtheorem{notation}{notation}
\newcommand{\bnt}{\begin{notation}}
\newcommand{\ent}{\end{notation}}

\newtheorem{prop}{Proposition}
\newcommand{\bprop}{\begin{prop}}
\newcommand{\eprop}{\end{prop}}

\newtheorem{lemma}{Lemma}
\newcommand{\blemma}{\begin{lemma}}
\newcommand{\elemma}{\end{lemma}}

\newtheorem{rmk}{Remark}
\newcommand{\brmk}{\begin{rmk}}
\newcommand{\ermk}{\end{rmk}}

\newtheorem{corr}{Corollary}
\newcommand{\bcorr}{\begin{corr}}
\newcommand{\ecorr}{\end{corr}}

\usepackage{mathtools}

\begin{document}

\title{On Koopman Operator for Burgers' Equation}% Force line breaks with \\

\author{Mikhael Balabane$^{1,2}$ \and Miguel A Mendez$^3$ \and 
	Sara Najem$^3$ }

\date{
	$^1$Laboratoire Analyse, G\'eom\'etrie et Applications, Universit\'e Paris 13\\
    $^2$Center for Advanced Mathematical Sciences, American University of Beirut\\
	$^3$von Karman Institute for Fluid Dynamics, EA Department, Sint Genesius Rode, Belgium\\
	$^4$Physics Department, American University of Beirut, Beirut 1107 2020, Lebanon
	%	\today
}


\maketitle

\begin{abstract}
 We consider the flow of Burgers' equation on an open set of (small) functions in $L^2([0,1])$. We \textcolor{black}{derive explicitly the Koopman decomposition of the Burgers' flow}. \textcolor{black}{We identify the frequencies and the coefficients of this decomposition as eigenvalues and eigenfunctionals of the Koopman operator. We prove the convergence of the Koopman decomposition for $t>0$ for small Cauchy data, and up to $t=0$ for regular Cauchy data. 
	The convergence up to $t=0$} leads to a `completeness' property for the basis of Koopman modes. 
We construct all modes and eigenfunctionals, including the eigenspaces involved in geometric multiplicity. This goes beyond the summation formulas provided by \cite{Page2018}, where only one term per eigenvalue was given. 
A numeric illustration of the Koopman decomposition is given and the Koopman eigenvalues compares to the eigenvalues of a Dynamic Mode Decomposition (DMD).

\end{abstract}
%\tableofcontents

\section{\label{sec:level1}Introduction}

The Koopman operator is a linear operator defined by Koopman in \cite{Koopman1931} to `linearize' nonlinear flows. This tool linearly evolves a set of observables (functionals) of the system state and allows for defining a spectrum for nonlinear flows. Such a spectrum is often considered as the generalization of `normal modes' for linear systems to `Koopman modes' for nonlinear systems \citep{Mezic2005,Mezic2013}. 

The decomposition of data on the Koopman modes is commonly referred to as Koopman Decomposition and offers promising applications in reduced-order modeling, feature extraction, and control \citep{Budisic2012,KutzDMD}. The Koopman framework is also invoked as a background for the Dynamic Mode Decomposition (DMD \cite{SCHMID2010}), a decomposition which finds the best linear dynamical system approximating snapshots of a flow (best in the least square sense). The link between DMD and Koopman decomposition was first revealed by \cite{ROWLEY2009} and is described by various authors (e.g., \textcolor{black}{\cite{Brunton2016,Rowley2017}}) on simple nonlinear ODEs. A systematic analysis of more complex problems is limited by the difficulties in deriving Koopman eigenfunctions for PDEs. The DMD is thus often considered as a data-driven method to compute Koopman modes under the assumption that the chosen observables and the available dataset are sufficiently `rich', in a sense defined by \cite{Tu2014}. Kernel methods to heuristically augment the richness of available data and observables have been proposed by \cite{Williams2015,Williams2015a} and reviewed by \cite{Kutz2018} in the framework of support vector machines (SVM).

The recent contribution by \cite{Page2018} has opened new avenues for  Koopman analysis.
These authors \textcolor{black}{proposed a procedure to derive} the Koopman decomposition for a nonlinear PDE, namely the Burgers' equation. \textcolor{black}{This procedure encounters a geometric multiplicity problem which is not solved. However, it} enabled the analysis of the (non-trivial) requirements for their data-driven identification via DMD. The main difficulties in identifying Koopman modes via DMD for problems with multiple invariant solutions are further discussed by \cite{Page2019}. More examples of Koopman analysis for PDEs are given by \cite{Nakao2020}.

In this work, we compute explicitly the set of eigen-observables of the Koopman operator for the \textcolor{black}{Burgers'} equation. We construct all modes and eigenfunctionals, including the eigenspaces involved
in geometric multiplicity. This goes beyond the summation formulas provided by
\cite{Page2018}, where only one term per eigenvalue was given. We prove convergence of the Koopman decomposition, and the completeness of the Koopman modes. %Moreover, we prove the decomposition of all Dirac functionals on the Koopman eigen-functionals. This is a way to prove ``completeness'' of the Koopman eigen-functionals when taking Dirac functionals as observables.

% The result can be extended to a Schr\"odinger equation.
% \textcolor{black}{ The same decomposition applies to a nonlinear Schr\"odinger equation, where the explicit computation of the Kolmogorov duplication of frequencies can be shown}.

The problem set and the relevant definitions are given in Section \ref{SecII}. \textcolor{black}{Section \ref{SecIII} recalls the Cole-Hopf transform, the methodology firstly proposed by Page and Kerswell \cite{Page2018}, and the geometric multiplicity problem they encounter. In section \ref{SecIV}, we propose a procedure that computes as many independent Koopman modes as needed by the encountered geometric multiplicity. We thus give an explicit Koopman decomposition. We link the coefficients of the Koopman decomposition to the eigenfunctionals of the Koopman operator. This leads to a `completeness' property of the Koopman eigenfunctionals. We give a numerical example of Koopman approximations in section \ref{SecVI} and compare its eigenvalues to those of a DMD.  Estimates needed in the proofs are given in section \ref{SecV}, where the smallness of the Cauchy data required for the decomposition's convergence is quantified. Concluding remarks are given in section \ref{SecVII}.}

\section{Definitions and Scope}\label{SecII}

We consider the Burgers' equation on $[0,1]$:

\begin{equation}
\label{F1}
\partial_t u =F_{\mathcal{B}}(u):= -u\partial_xu+\partial_{xx}u\quad u(t,0)=u(t,1)=0\quad u(0,x)=u_0(x)\,.
\end{equation}

Any Burgers' equation, with $x \in [0,L]$ and viscosity $\nu$, can be reduced to this form with a proper re-scaling of $x$, $u$, and $t$. Equation (\ref{F1}) translates to:

%  The Burgers flow $\Phi^t_{\mathcal{B}}(u_0)$ is defined on the set 

\begin{equation}
\label{flow}
u(t,x)=\Phi^t_{\mathcal{B}}(u_0)=u_0(x)+\int^{t}_{0} F_{\mathcal{B}} (u)(s)ds\,,
\end{equation}

\textcolor{black}{We define $\Omega_B$ as a subset of $L^2(]0,1[)$ invariant under $\Phi_{\mathcal{B}}^t$ and where a Koopman decomposition converges for $t>0$. The precise definition of $\Omega_B$ is given in Section \ref{SecIV} in formula (\ref{Omega_0})}. \noindent We consider the set $\mathcal {O}_{\mathcal{B}}$ of continuous observables (functionals) on $\Omega_{\mathcal{B}}$, i.e. the set of continuous maps $\phi$ from ${\Omega_{\mathcal{B}}}$ to $\mathbb{R}$. For any $t>0$, the Koopman operator $\mathcal{K}^t_{\mathcal{B}}$ of the Burgers' flow is the map from $\mathcal {O}_{\mathcal{B}}$ to itself defined by

\beq
\forall \phi \in \mathcal {O}_{\mathcal{B}}, \, \forall u_0 \in {\Omega_{\mathcal{B}}},\, (\mathcal{K}^t_{\mathcal{B}}(\phi))(u_0)=\phi(\Phi^t_{\mathcal{B}}(u_0))\label{F2}\,.
\eeq

\noindent $(\mathcal{K}_{B}^t)_{t\geq 0}$ is a flow on $\mathcal{O}_{B}$. $\mathcal{K}^t_{\mathcal{B}}$ is a linear map, and fulfills the multiplicative property

\beq \label{Fmult}\textcolor{black}{\forall \phi, \phi' \in \mathcal{O}_{\mathcal{B}},} \quad \forall u_0 \in {\Omega_{\mathcal{B}}}, \quad  (\mathcal{K}^t_{\mathcal{B}}(\phi\phi'))(u_0)= (\mathcal{K}^t_{\mathcal{B}}(\phi))(u_0)(\mathcal{K}^t_{\mathcal{B}}(\phi'))(u_0)\,.\eeq

The eigen-observables (eigen-functionals) $\varphi_\nu(u_0)$ \textcolor{black}{of the Koopman operator} are observables whose evolution by the Koopman flow is:

\beq
\label{eig}
(\mathcal{K}^t_{\mathcal{B}}(\varphi_\nu))(u_0)=\varphi_\nu(u_0) e^{\lambda_{\nu} t}\,,
\eeq
\textcolor{black}{
	\noindent with $\lambda_\nu$ the associated eigenvalue and $e^{\lambda_{\nu} t}$ the corresponding temporal evolution. Combining \eqref{Fmult} and \eqref{eig} shows that the product of two eigen-observables $\varphi_\nu$ and $\varphi_{\nu'}$, with eigenvalues $\lambda_\nu$ and $\lambda_{\nu'}$, is also an eigen-observable $\varphi_{\nu''}=\varphi_\nu\varphi_{\nu'}$ with eigenvalue $\lambda_{\nu''}=\lambda_\nu+\lambda_{\nu'}$.}
\medskip

The usual assumption in the Koopman decomposition is that the set of eigen-observables is sufficiently large to represent \emph{any} observable. Therefore, the value of \textit{any} observable $\phi$ at any function $f$ can be written as

%% MIKHAEL %%

\begin{equation}
\label{phi}
\phi(f)=\sum_{\nu } a_{\nu}\textcolor{black}{(\phi)}\varphi_{\nu}(f)  \,.
\end{equation}
\textcolor{black}{
	Combining \eqref{phi} with \eqref{F2} and \eqref{eig}, the observable $\phi$ evolves along a Burgers \textcolor{black}{orbit} through the Koopman flow: }

\begin{equation}
\label{Koopman_DE2}
\phi(\Phi^t_{\mathcal{B}}(u_0))=(K^t_{\mathcal{B}} \phi) (u_0)=\sum_{\nu } a_{\nu}(\phi) \varphi_{\nu}(\Phi_{\mathcal{B}}^t (u_0))=\sum_{\nu } a_{\nu}(\phi) (K^t_{\mathcal{B}}\,\varphi_{\nu})(u_0) =\sum_{\nu } a_{\nu}(\phi)\varphi_{\nu}(u_0) e^{\lambda_v t} \,.
\end{equation}

Equation \eqref{Koopman_DE2} is known as Koopman decomposition with respect to the observables $\phi$ and $a_\nu(\phi)$ are the associated Koopman modes. We consider the Dirac observables $\phi:=\delta_x$, where $\delta_x(u_0)=u_0(x)$: this observable maps any state variable $u$ to its value at location $x$. Equation \eqref{Koopman_DE2} becomes:

\begin{equation}
\label{Koopman_DE}
u(t,x)=\Phi^t_{\mathcal{B}}(u_0) =\sum_{\nu } e^{\lambda_v t} \varphi_{\nu}(u_0)  a_{\nu}(x)\,.
\end{equation}

While this decomposition looks linear, the nonlinearity of Burgers equation shows in the nonlinearity of the coefficients $\varphi_\nu(u_0)$.

We derive the Koopman  decomposition \eqref{Koopman_DE} for Burgers equation (\ref{F1}) for $u_0 \in \Omega_{\mathcal{B}}$. We prove its convergence in section \ref{SecIV}. We give explicit formulas for the Koopman modes $a_{\nu}(x)$, the eigenvalues $\lambda_{\nu}$ and the coefficients $\varphi_{\nu}(u_0)$ in (\ref{Koopman_DE}). We prove that $(e^{\lambda_\nu t},\varphi_\nu)$ are the eigen-elements of the Koopman operator associated with the Burgers equation. It is important to stress that the convergence of \eqref{Koopman_DE} requires smallness of the Cauchy data as shown in section \ref{SecIV}, hence the restriction of the Burgers flow to $\Omega_{\mathcal{B}}$. 
\section{Problem Statement via Cole-Hopf Transform}\label{SecIII}

\textcolor{black}{The Burgers equation \eqref{F1} is one of the few examples of nonlinear PDEs amenable to Koopman Analysis, thanks to the linearizing Cole-Hopf transformations. This was noticed by Kutz \emph{et al} \cite{Kutz2018} and exploited by Page and Kerswell \cite{Page2018}. }

\textcolor{black}{The Cole-Hopf transforms are defined as:}

\textcolor{black}{
	\beq
	u:=C(v)=-2{\partial_xv\over v}\quad {\rm and} \quad v:=H(u)={e^{-{1\over 2}\int_0^xu(s)ds}\over \int_0^1e^{-{1\over 2}\int_0^xu(s)ds}dx}\,.
	\label{colehopf}
	\eeq}
$H$ is defined for all functions in $L^2([0,1])$,  fulfills $\int_0^1H(u)(s)ds=1$, and $H(u) > 0$. We restrict $C$ to the set of functions $v$ having first weak derivative in $L^2([0,1])$, strictly positive, and with $\int_0^1v(s)ds=1$. This makes $H$ and $C$ inverse transforms.  %, and bounded away from zero. 
% We assume $\int_0^1v(s)ds=1$ and bounded\-ness of $1/ v$ follows from smallness of $\Vert \partial_xv\Vert_{L^2}$ as quantified in  the lemmi below. 

The Cole-Hopf transforms are  central tools because if $u(t,.)$ solves Burgers equation (\ref{F1}), then $v(t,.)=H(u(t,.))$ solves the linear heat equation:

\textcolor{black}{\begin{equation}
	\partial_tv=\partial_{xx}v\quad  \partial_xv(t,0)=\partial_xv(t,1)=0\quad v(0,x)=H(u_0):=v_0 \label{F5}\,,
	\end{equation}
	and the  converse is true. }

Writing the Cauchy data in the Fourier basis with $e_m(x):=\sqrt{2}\cos(m\pi x)$ for $m\not= 0$ as: $$ v_0(x) = 1 + \sum_{m=1}^\infty  c_m(v_0)e_m(x) \quad { \rm with} \quad c_m(v_0)=  \int_0^1 v_0(s) e_m(s) ds\,,$$
the solution of the heat equation is:

%% Comment for Miguel: You should introduce the c_m in the previous equation

\begin{equation}\label{vfourier}
v(t,x)=\Phi^t_{\mathcal{C}}(v_0)= 1 + \sum_{m=1}^\infty e^{-m^2\pi^2t} c_m(v_0)e_m(x)\,.
\end{equation} %=\sum_{m=1}^\infty e^{-m^2\pi^2t} l_m(u_0)e_m(x)
If $\Phi^t_{\mathcal{C}}$ denotes the flow associated with the heat equation, then $\Omega_{\mathcal{C}}:=H(\Omega_{\mathcal{B}})$ is invariant under this flow and we have:

\begin{equation}
\textcolor{black}{
	H(\Phi^t_{\mathcal{B}}(u_0))=\Phi^t_{\mathcal{C}}( H(u_0))\quad {\rm and}\quad  C(\Phi^t_{\mathcal{C}}(v_0))=\Phi^t_{\mathcal{B}}(C(v_0))}\,.
\label{F6_M}
\end{equation}

\textcolor{black}{Using the composition law $(A\circ B)(z):= A(B(z))$, }equation \eqref{F6_M} can be written as

\begin{equation}
H\circ \Phi^t_{\mathcal{B}}=\Phi^t_{\mathcal{C}}\circ H\quad {\rm and}\quad  C\circ \Phi^t_{\mathcal{C}}=\Phi^t_{\mathcal{B}}\circ C\,.
\label{F6}
\end{equation}

These equations give the conjugacy properties linking the flows $\Phi^t_{\mathcal{B}}$ and $\Phi^t_{\mathcal{C}}$ via the Cole-Hopf transform \eqref{colehopf}. This intertwining of the Burgers flow with the linear heat flow gives the Koopman decomposition for Burgers' equation.

Page and Kersell's procedure in \cite{Page2018} consists in writing \eqref{colehopf} as $u v = -2\partial_x v$, plugging the Fourier decomposition \eqref{vfourier} of $v$ and a formal Koopman decomposition \eqref{Koopman_DE} of $u$, and then identifying the Koopman modes, amplitudes and eigenvalues by inspection. Except for modes with multiplicity one, this procedure encounters a problem. 

In the next sections, we propose an approach which solves the issue of geometric multiplicity. We explicitly give the Koopman decomposition and prove its convergence.

\section{An Explicit Koopman Decomposition}\label{SecIV}

\subsection{An Analytical Decomposition into Exponentials}\label{SubSecIV_I}

\noindent  For any $u_0\in \Omega_{\mathcal{B}}$, let  $v_0=H(u_0)$. \textcolor{black}{We note that the steady state solution $v=1$ is a sink for the heat equation, i.e.  $v(t,x)\rightarrow1$ for $t\rightarrow +\infty$ for all $v_0 \in \Omega_C$.} \textcolor{black}{We restrict our decomposition to the vicinity of this steady state and write $v_0=1+w_0$, with $w_0$ a function such that $\int_0^1w_0(x)dx=0$ and $\Vert \partial_x w_0 \Vert_{L^2}\leq {1\over 4}$. From Section \ref{lemmi} and the choice of $\Omega_B$, we get that $\Vert \partial_x w_0 \Vert_{L^2}\leq {1\over 4}$, hence  $\sup_x\vert w_0(x)\vert \leq 1/4$ (by property 1 in section \ref{lemmi})}.

\medskip

\noindent For $t>0$, \textcolor{black}{we} write $\Phi^t_{\mathcal{C}}(v_0)=v(t,.)=1+ w(t,.)$. $\partial_xw$ fulfills the heat equation \textcolor{black}{ with Dirichlet  boundary conditions. Energy decay due to diffusion gives} $\Vert \partial_x w\Vert_{L^2}\leq \Vert \partial_x w_0\Vert_{L^2}$ so $\sup_x\vert w(t,x)\vert \leq {1\over 4}$. One can therefore apply the Cole-Hopf transforms, \textcolor{black}{and have the following asymptotic expression for the solution of the Burgers's equation:}
\beq \label{Fdev}\Phi^t_{\mathcal{B}}(u_0)=u(t,.)=C(v(t,.))=-2{ \partial_x w\over {1+ w}}=-2 \partial_x w\sum_{q=0}^\infty (-1)^q w^q\,.
\eeq

\textcolor{black}{Using the expression for $w$ given by formula (\ref{vfourier}) we compute} $w^q$, the product of $q$ identical sums, by taking  term \textcolor{black}{of rank} $n_1$ in the first sum, term \textcolor{black}{of rank} $n_2$ in the second sum, \textcolor{black}{and so on up to $m$-th} sum. We get:
$$w^q=\sum_{(n_1,..,n_q)\in {(\textcolor{black}{N\backslash  \{0\}})}^q }e^{-\sum_{k=1}^qn_k^2 \pi^2t}\prod_{k=1}^{q}c_{n_k}(v_0)\prod_{k=1}^{q}e_{n_k}(x)\,.
$$ \noindent \textcolor{black}{On the other hand:}

$$\partial_xw(t,x)=\sum_{n_0=1}^\infty  e^{-n_0^2\pi^2t} c_{n_0}(v_0)\partial_xe_{n_0}(x)\,.
$$
\noindent
Replacing the above expressions in formula (\ref{Fdev}) gives
\begin{equation}
\begin{split}
\Phi^t_{\mathcal{B}}(u_0)(x)=-2
\sum_{n_0=1}^\infty  e^{-n_0^2\pi^2t} c_{n_0}(v_0)\partial_xe_{n_0}(x)
\\-2\sum_{q=1}^\infty \sum_{n_0=1}^\infty \sum_{\stackrel{n_1,..,n_q}{\in {(\textcolor{black}{N \backslash \{0\}} )}^q}}
(-1)^{q}e^{-\pi^2 t\sum_{k=0}^qn_k^2} \prod_{k=0}^{q}c_{n_k}(v_0)\,\,\partial_xe_{n_0}(x)\prod_{k=1}^{q}e_{n_k}(x)\,.
\end{split}
\end{equation}

\textcolor{black}{This decomposition is more easily handled by taking the following set of indices:}

\begin{equation}
\mathcal{A} = \{\nu = (n_0,n_1,...n_{\alpha(\nu)}); \alpha(\nu) \in \mathbb{N}, n_0 \in \mathbb{N}, n_i \in \mathbb{N}\backslash \{0\} \,\,{\rm for}\,\,i=1,..,\alpha(\nu)\,.
\end{equation}\label{indices}
Then, the above formula for $\Phi^t_{\mathcal{B}}(u_0)$ looks familiar if one introduces: 
\begin{equation}
\label{lambdas}
\lambda_\nu=-\pi^2\sum_{k=0}^{\alpha(\nu)}n_k^2\,,
\end{equation}
\textcolor{black}{
	\begin{equation}
	a_{\nu}(x)=(-1)^{\alpha(\nu)}2^{\alpha(\nu)+3\over 2}
	n_0\pi \sin{(n_0\pi x)}\prod_{k=1}^{\alpha(\nu)}{\cos{(n_k\pi x)} }
	\label{F_a_nu}\,,
	\end{equation} }
\textcolor{black}{
	\begin{equation}
	\varphi_{\nu}(u_0)=\prod_{k=0}^{\alpha(\nu)}l_{n_k}(u_0)\quad {\rm with}\quad l_n(u_0)=c_n(H(u_0))\,.
	\label{Phi_nu}
	\end{equation} }

\noindent \textcolor{black}{With these notations, we have derived the following Koopman decomposition:}

% $$\Phi^t_{\mathcal{B}}(u_0)(x)=\sum_{\nu \in \mathcal{A}}e^{-\lambda_\nu t}\varphi_\nu(u_0)a_\nu(x)
% $$
\textcolor{black}{
	\begin{equation}
	\forall u_0\in \Omega_{\mathcal{B}},\quad \forall t>0 \quad \forall x\in [0,1]\quad \Phi^t_{\mathcal{B}}(u_0)(x)=u(t,x)=\sum_{\nu\in  \mathcal{A}}  e^{\lambda_\nu t}\,\varphi_{\nu}(u_0)a_{\nu}(x)\,.
	\label{F3}\end{equation} }

\noindent \textcolor{black}{
	Note that the length of $\nu$ is $\alpha(\nu)+1$ and in formula (\ref{F_a_nu}) the last product should be taken as $1$ for $\alpha(\nu)=0$.} 

\textcolor{black}{The proper set of initial conditions $u_0$ granting convergence of this series is $\Omega_\mathcal{B}$, that we now define:}
\textcolor{black}{
	\begin{equation}
	\Omega_\mathcal{B}=\{u_0\in L^2([0,1]); 2e^{\Vert u_0\Vert_{L^2}}\Vert u_0\Vert_{L^2}<1\}\,.
	\label{Omega_0}
	\end{equation} }
\textcolor{black}{Section \ref{proof2} shows that the convergence of \eqref{F3} is uniform and in absolute values.}
% Convergence of the series \eqref{F3} is uniform, and in absolute values: this is proved in section \ref{proof2}.

% \textcolor{black}{It is worth noting that in formula (\ref{F3}), in order to have the least number of harmonics in the dynamics, one should have, for $\alpha_{\nu} = 0$, only one non-zero coefficients. Let $u_0$ be such that $l_n(u_0)= 0$ for $n_0 \neq p_0$. Then the only non-zero coefficient in the decomposition are those whose indices contain integers equal to $p_0$. The decomposition then writes: } 
% \textcolor{black}{
% \begin{equation} \label{F3bisbis}    u(t,x)=\sum_{k=0}^{\infty} (-1)^k 2^{\frac{k+3}{2}} p_0 \pi e^{-k p_0^2 \pi^2 t}(l_{p_0}(u_0))^k\sin{(p_0 \pi x) } (\cos{(p_0 \pi x)})^k 
% \end{equation}}

% \textcolor{black}{This shows that all harmonics of the fundamental frequency do appear, with amplitudes that behave as a power. It seems this is the first instance where Kolmogorov duplication of frequency can be computed explicitly.  }

\subsection{The Koopman Eigenfunctionals for Burgers' Equation}\label{SubSecIV_II}

\textcolor{black}{We define $\mathcal{O}_{\mathcal{C}}$ as the set of observables on $\Omega_\mathcal{C}$. Let $\mathcal{K}_{\mathcal{C}}^t$ be the Koopman operator of the heat equation \eqref{vfourier}}.
The conjugacy formula (\ref{F6}) leads to the following adjoint identities linking $\mathcal{K}^t_{\mathcal{B}}$ and $\mathcal{K}^t_{\mathcal{C}}$:
% Here you should say something more about these conjugacy

\textcolor{black}{\begin{subequations} 
		\label{F7}
		\begin{equation}
		\label{F7a}
		\forall \phi \in \mathcal{O}_{\mathcal{B}}\,,\,\, (\mathcal{K}^t_{\mathcal{B}}\phi)\circ C=\mathcal{K}^t_{\mathcal{C}}(\phi\circ C)\,,
		\end{equation}
		\begin{equation}
		\label{F7b}
		\forall \psi\in \mathcal{O_{\mathcal{C}}}\,,\,\, (\mathcal{K}^t_{\mathcal{C}}\psi)\circ H=\mathcal{K}^t_{\mathcal{B}}(\psi\circ H)\,.
		\end{equation}
\end{subequations} }
Using previous notations, formula \eqref{F7a} is true because of the following:
% \textcolor{black}{$$ (\mathcal{K}^t_{\mathcal{B}} \phi)(C(v_0))  = \phi (\Phi_{\mathcal{B}}^t(C(v_0))) = \phi((\Phi^t_{\mathcal{B}} \circ C)(v_0)) =\phi((C \circ \Phi^t_{\mathcal{C}})(v_0) )=(\phi \circ C)(\Phi^t_{\mathcal{C}}(v_0))=(\mathcal{K}^t_{\mathcal{C}}(\phi \circ C))(v_0)$$}

$$(\mathcal{K}^t_{\mathcal{B}} \phi)(C(v_0)) = (\mathcal{K}^t_{\mathcal{B}} \phi)(u_0) = \phi(\Phi^t_{\mathcal{B}}u_0) = \phi(u(t,.)) =  \phi(C (v(t,.))) = (\phi \circ C)(v(t,.)) = (\phi \circ C) ( \Phi^t_{\mathcal{C}}(v_0)) = \mathcal{K}^t_{\mathcal{C}}(  \phi \circ C) (v_0)\,.$$

This also applies to \eqref{F7b}.

We now prove that the functions $\varphi_\nu(u_0)$ and the exponents $\lambda_\nu$ of decomposition (\ref{F3}) are eigenfunctionals and the associated eigenvalues of the Koopman operator in formula \eqref{F2}, i.e. they satisfy property (\ref{eig}).

%\bprop \label{prop1} 
The proof goes as follows:
\textcolor{black}{the linearity of $c_n$ and $\Phi^t_{\mathcal{C}}$ in \eqref{vfourier} gives }

\begin{equation}
(\mathcal{K}^t_{\mathcal{C}}(c_n))(v_0) = c_n(\Phi^t_{\mathcal{C} }(v_0)) = c_n \biggl(1 + \sum_{m=1}^\infty e^{-m^2\pi^2t} c_m(v_0)e_m(x)\biggr)=c_n(v_0)e^{-n^2\pi^2 t} \,,   
\end{equation}

% We take $v_0=H(u_0)$ and use \eqref{F7b} to get that all observables $l_n$ in \eqref{Phi_nu} are eigenfunctionals of $\mathcal{K}^t_{\mathcal{B}}$: }
\textcolor{black}{so, $\mathcal{K}^t_{\mathcal{C}}(c_n) = e^{-n^2 \pi^2 t} c_n$,  and by \eqref{F7b} $\mathcal{K}^t_{\mathcal{B}}(c_n \circ H )= e^{-n^2 \pi^2 t }(c_n \circ H)$. Then:}

\beq \label{F9}
\forall n\in N \quad \mathcal{K}_{\mathcal{B}}^t(l_n)=e^{-n^2\pi^2t}l_n\,.
\eeq
% \eprop

%\bcorr \label{corr1}
The multiplicative property (\ref{Fmult}) gives, for all $\nu=(n_0,..,,n_{\alpha(\nu)}) \in \mathcal{A}$, 

\beq \label{Ffi}
\mathcal{K}^t_{\mathcal{B}}(\prod_{k=0}^{\alpha(\nu)} l_{n_k}) = e^{\lambda_{\nu} t} \prod_{k=0}^{\alpha(\nu)} l_{n_k}. \eeq
%\ecorr
\noindent \textcolor{black}{This proves that  $\lambda_{\nu}$ is an eigenvalue of ${\mathcal{K}}^t_{\mathcal{B}}$} with \textcolor{black}{eigen-functional} $\varphi_\nu$.
\textcolor{black}{This identifies coefficients and exponentials in formula (\ref{F3}) as spectral elements of $\mathcal{K}^t_{\mathcal{B}}$: the coefficients are the value taken by the eigenfunctionals at the Cauchy data.}
\textcolor{black}{The $a_{\nu}$ will be identified in the next section}. 

Notice that the definition of these eigenfunctionals needs no assumption on the smallness of the state variables.

\subsection{On completeness}\label{SubSecIV_III}

\textcolor{black}{We call weak-completeness of these eigenfunctionals the property that any Dirac functional $\delta_x$ can be decomposed on these eigenfunctionals.} \textcolor{black}{To prove this weak-completeness property},  it is worth examining under what assumption
the convergence of formula (\ref{F3}) is valid at $t=0$. \textcolor{black}{A suitable assumption is a regularity assumption on initial conditions: square integrability of the first derivative of the Cauchy data.  That leads to define the following set of initial conditions: }

$$\omega_{\mathcal{B}}=\biggl \{u_0\in \Omega_{\mathcal{B}};\quad u_0(0)=u_0(1)=0;\quad \int^1_0\vert \partial_xu_0\vert^2<\infty\biggr\}.$$

We prove in section \ref{proof3} that for $u_0\in\omega_{\mathcal{B}}$ the convergence of formula (\ref{F3}) is uniform, for $t\geq 0$ and $x\in [0,1]$. This implies:

%\bcorr  
\beq \forall u_0\in \omega_{\mathcal{B}},\quad \forall x\in [0,1], \quad u_0(x)=\sum_{\nu \in  \mathcal{A}} \varphi_{\nu }(u_0) a_{\nu}(x).\label{F4}\eeq
%\ecorr

If one takes for observables the values at specific locations $x$, denoted by the Dirac notation $\delta_x$, formula (\ref{F4}) can be written as: 
% \noindent By noticing that if $u_0\in \omega_B$ then $u_0$ is a continuous function, formula (\ref{F4}) can be re-written.
% \noindent One can better state formula (\ref{F4}) by noticing that if $u_0$ has a square integrable derivative then $u_0$ is a continuous function. 
% Hence for any $x\in [0,1]$, the Dirac observable $\delta_x$, which maps a state variable $u$ to its value at location $x$, is a continuous observable on $\omega_{\mathcal{B}}$. Then, Corollary \ref{corr2} can be written in $\mathcal{O}_{\omega_B}$ as:

%Then the previous corollary can be re-written as:

% \bcorr \label{corr3} 

\beq 
\forall x\in [0,1]\quad \delta_{x}=\sum_{\nu \in  \mathcal{A}} a_{\nu}(x)\varphi_{\nu }.
\label{F4bis}\eeq
% ecorr

%\noindent Here, for any given $x\in [0,1]$, the sequence $(a_{\nu}(x))_{\nu \in \mathcal{A}}$ is the sequence of coefficients in the decomposition of $\delta_x$ on the set of Koopman eigen-observables. This can be thought of as `completeness' of the eigen-observables $(\varphi_\nu)_{\nu \in \mathcal{A}}$ although lack of a topology on the set ${\mathcal{O}}_{\omega_B}$  of observables makes completeness irrelevant.

% this to replace the last paragraph of 'main results'  beginning at 'here for any given $x\in [0,1]$', so the word 'irrelevant that bothers miguel desappears, and the energy formula is given.

% \noindent Here, for any given $x\in [0,1]$, the sequence $(a_{\nu}(x))_{\nu \in \mathcal{A}}$ is the sequence of coefficients in the decomposition of $\delta_x$ on the set of Koopman eigen-observables. Corollary \ref{corr3} proves \eqref{phi} for $\phi:=\delta_x$. This can be thought of as `completeness' of the eigen-observables $(\varphi_\nu)_{\nu \in \mathcal{A}}$, although  a set  of observables can be endowed only with the weak topology (convergence of values of observables at any given state variable), and is not normed. 
\textcolor{black}{Comparing this with \eqref{phi}, we identify the Koopman modes $a_{\nu}(x)$ as the coefficients of the decomposition of $\delta_x$ on eigenfunctionals of the Koopman operator.}

\medskip 

To illustrate this weak-completeness, one can show that formula \eqref{F4bis} implies the decomposition of the kinetic energy observable:

$$E(u)=  \int_0^1u^2(s)ds=\int_0^1(\delta_s(u))^2ds.$$ 

For $\nu=(n_0,\cdots n_m)\in \mathcal{A}$ and $\nu'=(n'_{0},\cdots n'_{m})\in \mathcal{A}$, we define the concatenation of these indices as

$$
c (\nu,\nu ')=\begin{dcases}
(n_0, \dots, n_m, n'_{0},\cdots n'_{m}) \quad \mbox {if} \quad n'_0\neq 0\\
(n_0, \dots, n_m, n'_{1},\cdots n'_{m}) \quad \mbox {if} \quad n'_0= 0
\end{dcases}
$$

\textcolor{black}{The multiplicative property of  the Koopman eigen-observables, 
	$\varphi_{\nu} \varphi_{\nu'} =\varphi_{c(\nu,\nu')}$ and the uniform convergence in the $x$ variable at $t=0$ of formula (\ref{F4})  implies: }

$$E(u_0)= \textcolor{black}{ \sum_{\mathcal{A} \times \mathcal{A}} \varphi_{\nu}(u_0) \varphi_{\nu'}(u_0) \int_0^1 a_{\nu}(s) a_{\nu'}(s) ds }= \sum_\mathcal{A\times A}b_{c(\nu,\nu')}\varphi_{c(\nu,\nu')}(u_0) \quad {\rm with}\quad b_{c(\nu,\nu')}=\int_0^1a_\nu(s)a_{\nu'}(s)ds\,.
$$

\noindent   The Koopman decomposition \eqref{Koopman_DE} \textcolor{black}{of the kinetic energy is then, by \eqref{Koopman_DE2}: } 
\textcolor{black}{
	$$E(u(t,.))=\sum_\mathcal{A\times A}e^{\lambda_{c(\nu,\nu')}t}b_{c(\nu,\nu')}\varphi_{c(\nu,\nu')}(u_0).
	$$}

\section{Numerical Illustrations of Koopman Approximations}\label{SecVI}

\textcolor{black}{
	We present a numerical illustration of the Koopman decomposition \eqref{F3} for the Burger's flow \eqref{F1} on the interval $[0,1]$. We consider initial conditions with $c_m=0$ for $m>2$ in \eqref{vfourier}:}

\begin{equation}
\quad u_0(x)=C(v_0)\quad \mbox{with} \quad v_0(x)=1+{1\over2}\cos{\pi x}+{1\over4}\cos{2\pi x}.
\label{C1}
\end{equation}
\textcolor{black}{
	The Burger's flow is computed from \eqref{colehopf} as $u(t,x)=C(v(t,x))$ with $v(t,x)$ obtained from \eqref{vfourier}.
	The computations are performed on a uniform mesh $x_i=(i-1)\Delta x$ with $i\in[1,1024]$, computing all integrals via the trapezoidal method. }

\textcolor{black}{
	We consider the set of Koopman modes originating from the indices of length one, up to length six ($0\leq\alpha(\nu)\leq 5$). This includes all the sets $\nu=(n_0)$, $\nu=(n_0,n_1)$, up to $\nu=(n_0,n_1,n_2,n_3,n_4,n_5)$, with $n_k\in[1,2]$ for all $k\leq5$. We thus consider a total of $126$ Koopman modes with repetitions, out of which $30$ are independent. All indices with $n_k>2$ leads to zero amplitude according to \eqref{Phi_nu} for this particular initial conditions given in \eqref{C1}, where $c_m=0$ for $m>2$.}

\begin{figure}[h!]
	\centering
	\subfigure[\ ]{
		\includegraphics[trim=5 30 5 30,clip,width=8.5cm]{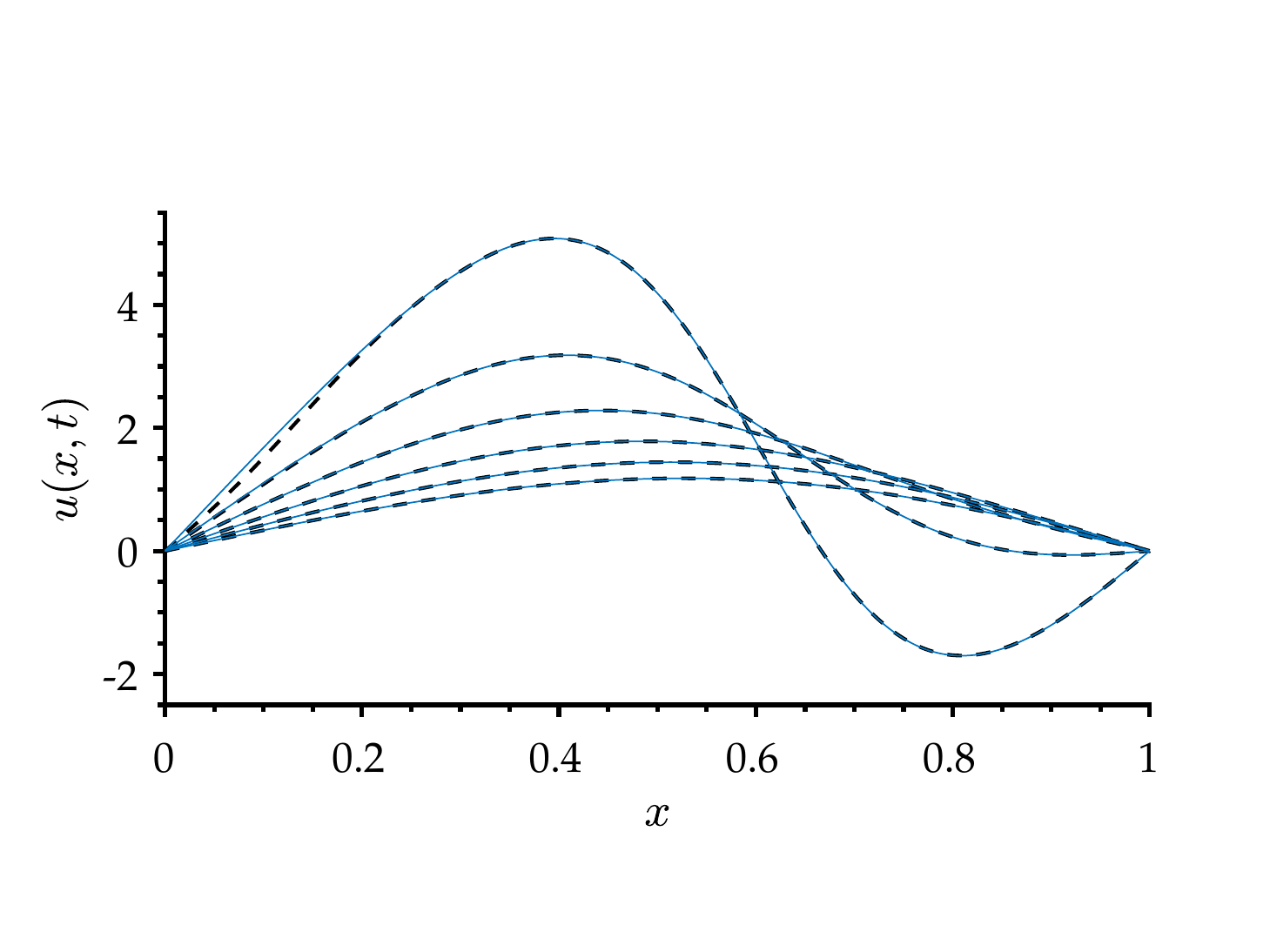}}
	\subfigure[\ ]{
		\includegraphics[trim=5 30 5 30,clip,width=8.5cm]{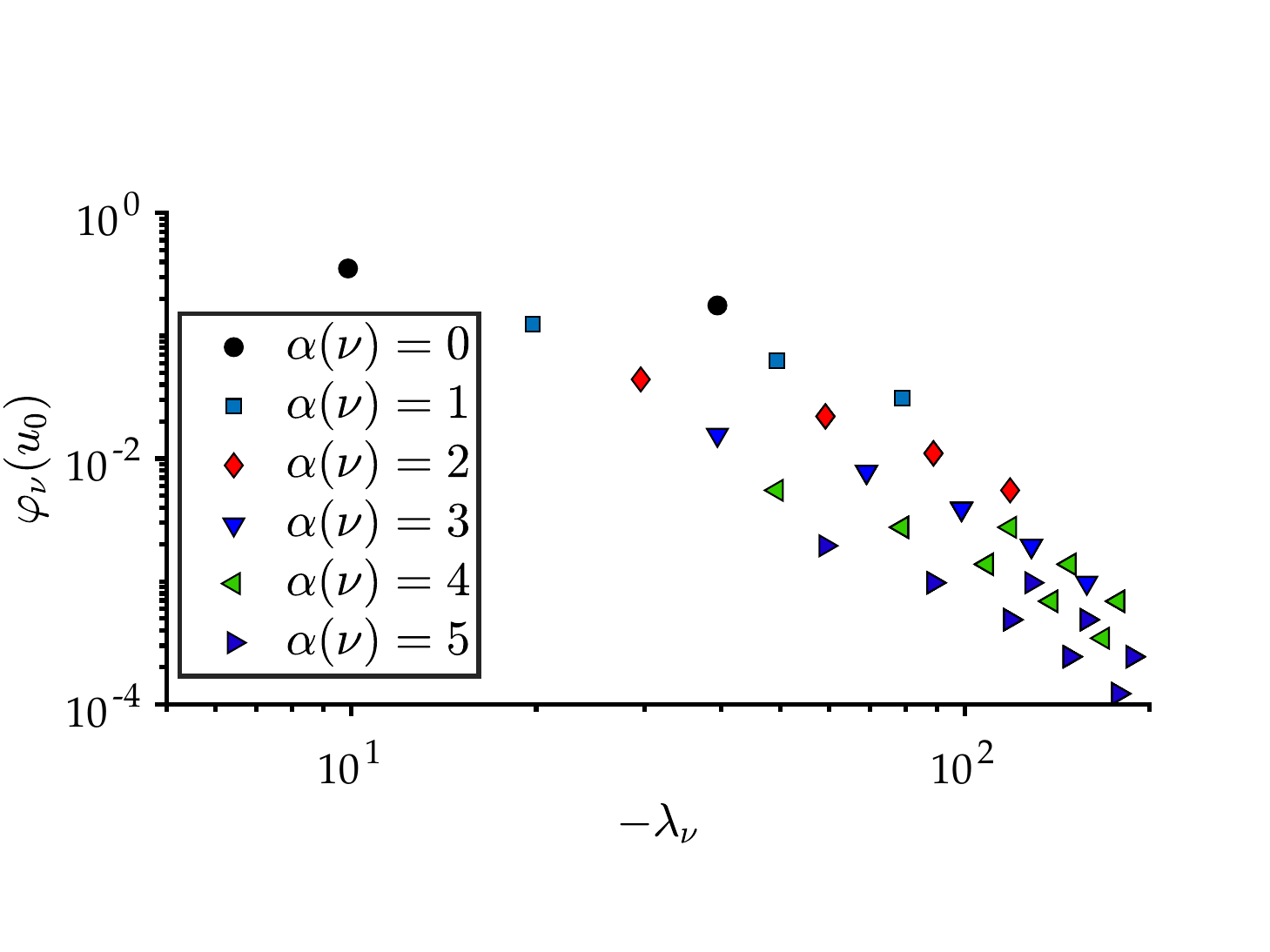}}
	\vspace{-2em}
	\caption{Koopman approximations of the Burgers' flow with initial condition in \eqref{C1}. Figure (a) shows the dynamics $u(x,t)$ (continuous blue line) along with the Koopman approximations (dashed lines). This initial condition does not satisfy the smallness condition needed in the proofs and this is why it does not converge up to $t=0$; Figure (b) plots the amplitudes of the Koopman modes $\varphi_\nu(u_0)$ as a function of $-\lambda_\nu$. }
	\label{Mig1}
\end{figure}

Figure \ref{Mig1}(a) shows the dynamics $u(x,t)$ (blue continuous lines) for $t=0,0.02, 0.04, 0.06, 0.14, 0.24$ along with the Koopman approximation (black dashed lines). While the completeness of the (infinite) Koopman basis has been proven, these results highlight the challenges in the convergence up to $t\rightarrow 0$: because $u_0$ does not fulfil the smallness property $u_0 \in \Omega_B$, convergence does not extend to $t=0$, as shown by the part of the Cauchy data close to $x=0$ (zoomed axis).

\textcolor{black}{
	Figure \ref{Mig1}(b) shows the Koopman coefficients $\varphi_{\nu}(u_0)$ as a function of the associated Koopman eigenvalue $\lambda_\nu$. The markers in the figures, labeled in the legend, recall the length of the set $\nu$ to which each mode corresponds. Notice that only 30 markers are visible. The multiplicity of the eigenvalues (hence the temporal evolution) depends on the set of indices (e.g, $\nu=(2)$ and $\nu=(1,1,1,1)$ lead to the same $\lambda_\nu$) and increases with $-\lambda_v$. The problem encountered by Page \& Kerswell \cite{Page2018} is that they identify the Koopman modes by combining their formulas (8) and (9) to get (10) and construct only one Koopman mode irrespective of whether its multiplicity is greater than one.}

\textcolor{black}{
	Because of the exponential evolution of the Koopman modes, the amplitudes $\varphi_\nu(u_0)$ weigh the contribution of each mode in the decomposition of the initial data. To analyze the relative contribution to the Burgers' flow $u(t,x)$, and to consider the role of the spatial basis $a_\nu (x)$, we introduce a measure of the mode relevance in the form of relative $L^2$ norm in space and time:}

\begin{equation}
\label{Eq_Importance}
\sigma_{\nu}=\frac{||e^{-\lambda_{\nu}t} \varphi_{\nu}(u_0)a_{\nu}(x) ||}{||u(t,x)||} \quad \mbox{with} \quad ||f||^2=\int^{t_2}_{t_1}\int_0^1 \bigl(f(x,t)\bigr)^2 dx \,dt\,.
\end{equation}

\begin{figure}[h!]
	\centering
	\subfigure[\ ]{
		\includegraphics[trim=5 30 5 30,clip,width=8.5cm]{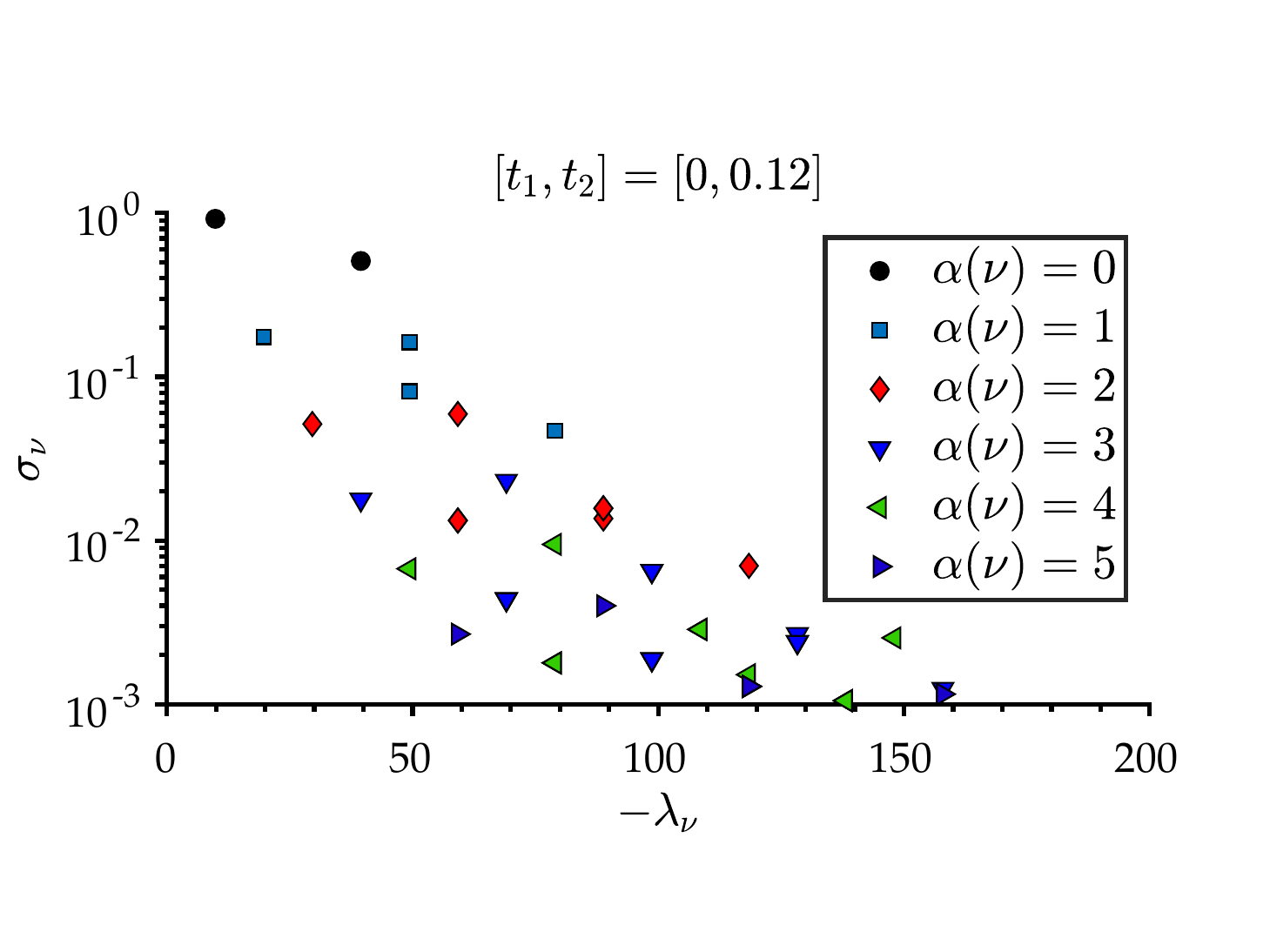}}
	\subfigure[\ ]{
		\includegraphics[trim=5 30 5 30,clip,width=8.5cm]{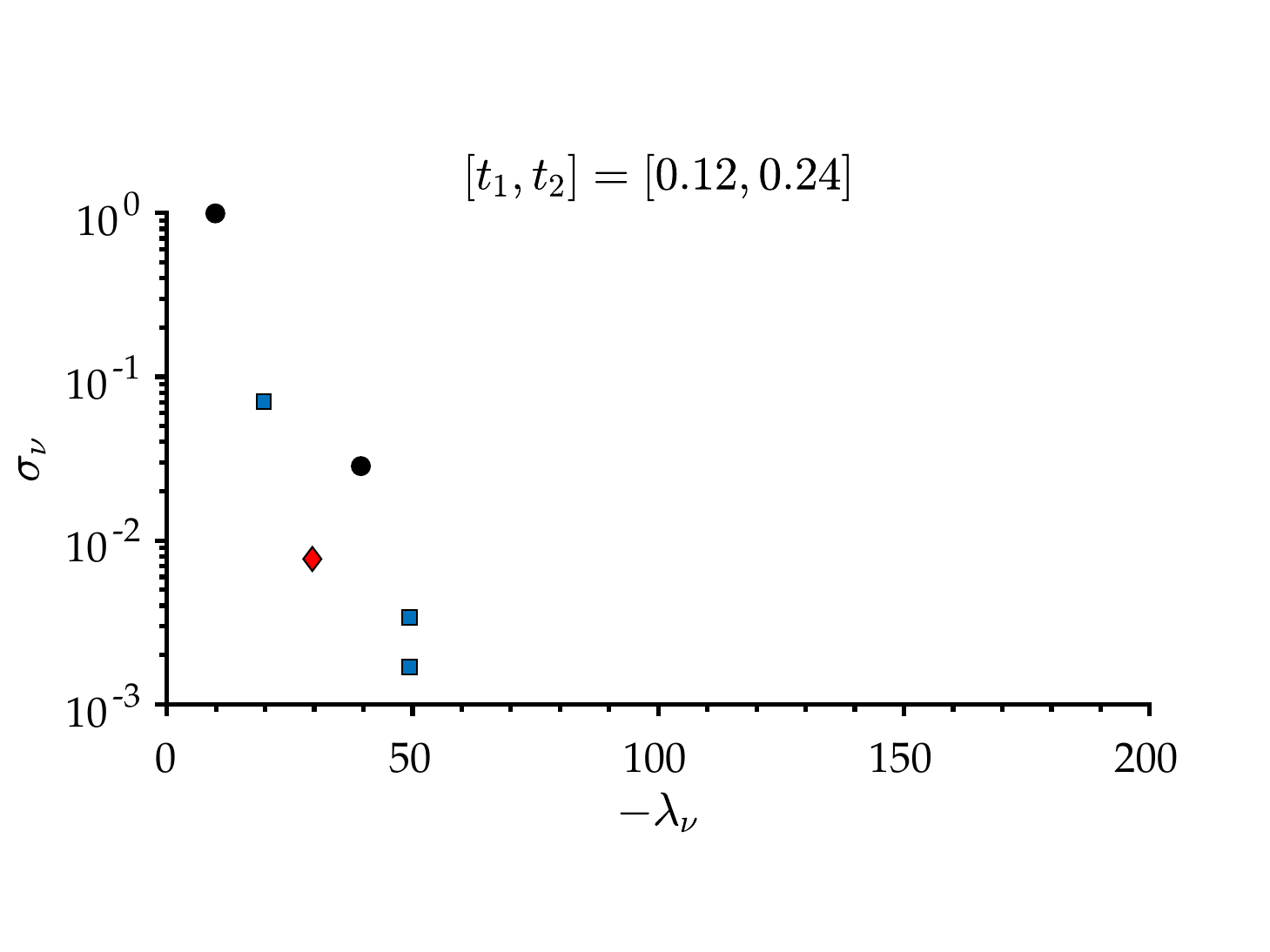}}
	\vspace{-1em}
	\caption{Contribution of each Koopman mode according to \eqref{Eq_Importance}, taking $[t_1,t_2]=[0,0.12]$ (a) and  $[t_1,t_2]=[0.12,0.24]$ (b).}
	\label{Mig2}
\end{figure}

\textcolor{black}{
	It is worth highlighting that this parameter only serves illustrative purposes, as the lack of orthogonality of the Koopman basis does not allow to recover the energy of the flow by summing the contribution of each mode.}
\textcolor{black}{
	The relative importance of the modes is analyzed for two time intervals, namely $[t_1,t_2]=[0,0.12]$ and $[t_1,t_2]=[0.12,0.24]$. The contribution of each mode according to \eqref{Eq_Importance} is shown in Figure \ref{Mig2}. Only modes with $\sigma_\nu>10^{-3}$ are considered.}

\textcolor{black}{
	A comparison of these figures shows how the relative importance of the Koopman modes changes in time, depending on the relative weight of advection and diffusion. Seven modes have $\sigma_\nu>0.05$ in the first time interval while only two have the same relevance in the second interval. The same plot considering $[t_1,t_2]=[0,24]$ is indistinguishable from Fig \ref{Mig2}(a) and it is thus not shown. This practically illustrates how the initial conditions lead to the definition of the decomposition. Moreover, this figure shows how adding the contribution of the spatial structure $a_\nu(x)$ results in the spreading of the plot in Figure \ref{Mig1}(a), as modes having the same amplitude ($\varphi(u_0)$) and temporal evolution ($\lambda_\nu$) might have different $a_\nu(x)$.} 

\begin{figure}[htbp]
	\centering
	\includegraphics[width=13cm]{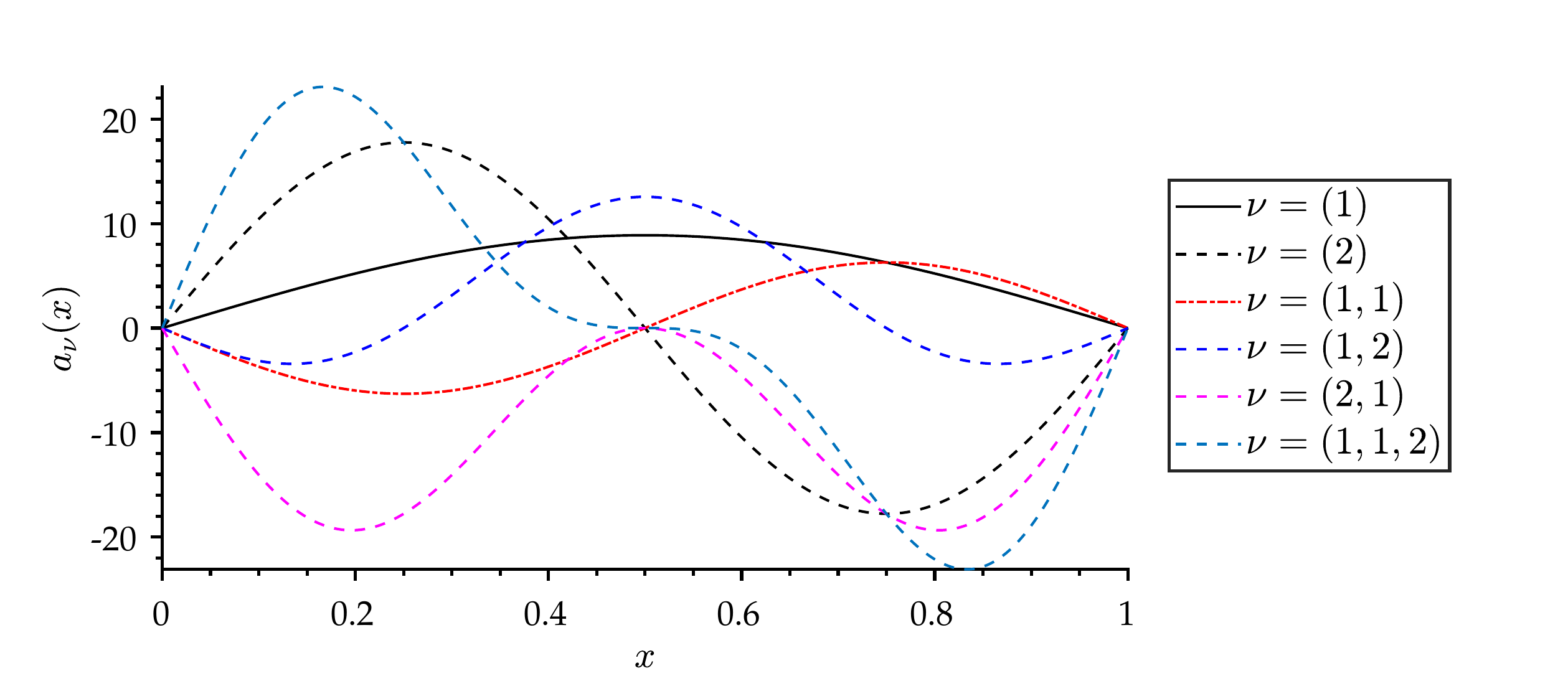}
	\caption{Five Koopman modes $a_\nu(x)$ of the Burgers' flow. Indices are indicated in the legend.}
	\label{a_nus}
\end{figure}

\textcolor{black}{
	The seven leading modes (with $\sigma_\nu>0.05$) taking the full time interval are shown in Figure \ref{a_nus}, with the legend recalling the associated set of indices $\nu$. Except for the first two, all the other modes display nonlinear interactions through the multiplicative property: for instance, the interaction of a mode with itself, that is self-interaction, gives rise to a squared mode.}

Finally, we conclude this illustrative section by analyzing the accuracy of a data-driven approximation of these Koopman modes using a Dynamic Mode Decomposition (DMD) (see \cite{SCHMID2010,ROWLEY2009,Tu2014}). We consider the SVD-based approach. A set of $n_t=101$ snapshots uniformly sampled in time discretization $t_k=(k-1)\Delta t$ with $\Delta t=0.002$ is used to construct the dataset matrix $\mathbf{X}\in\mathbb{R}^{1023\times 101}$.

\begin{figure}[!htbp]
\centering
\includegraphics[trim=5 30 5 30,clip,width=8.5cm]{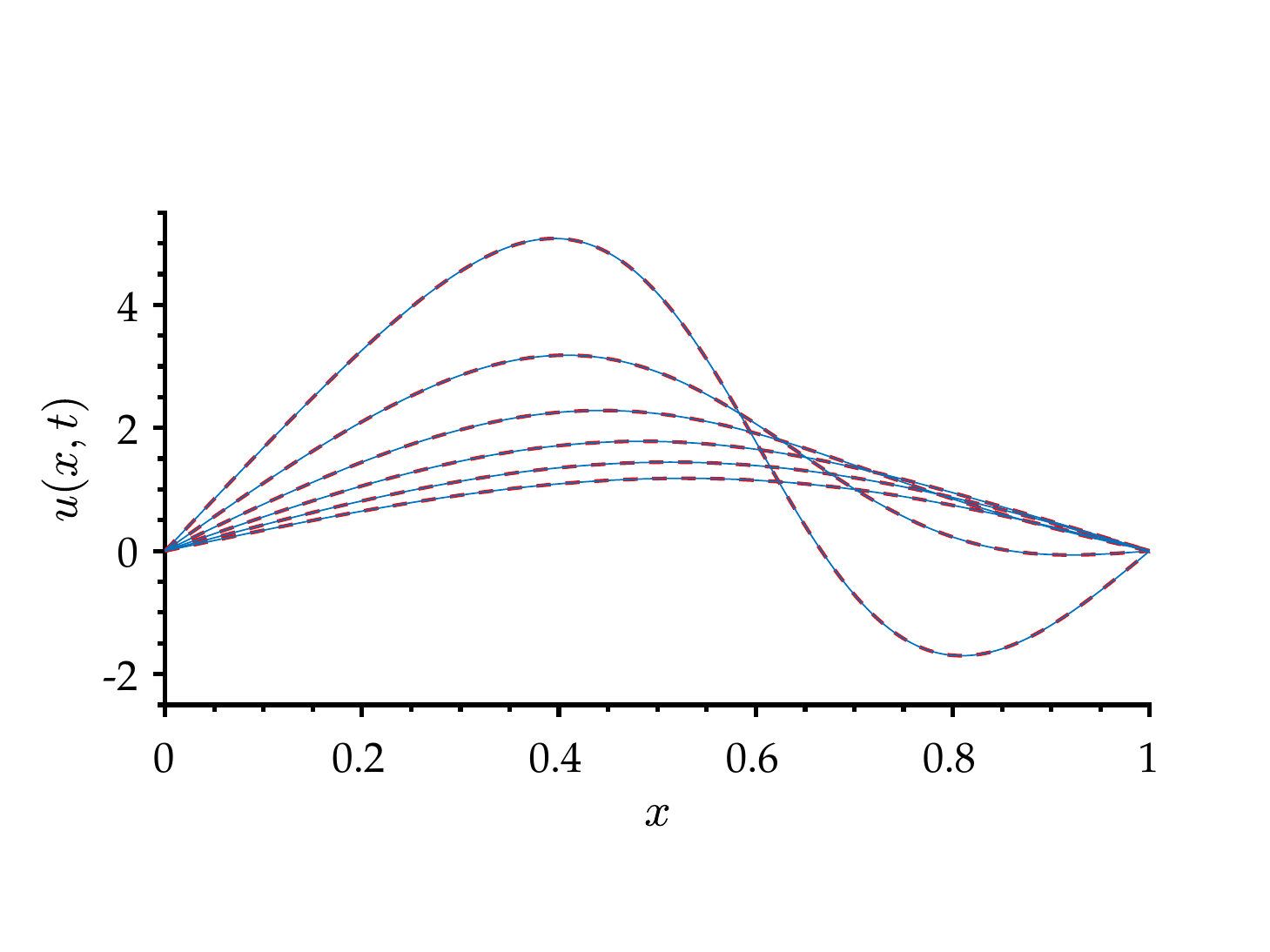}
\includegraphics[width=8cm]{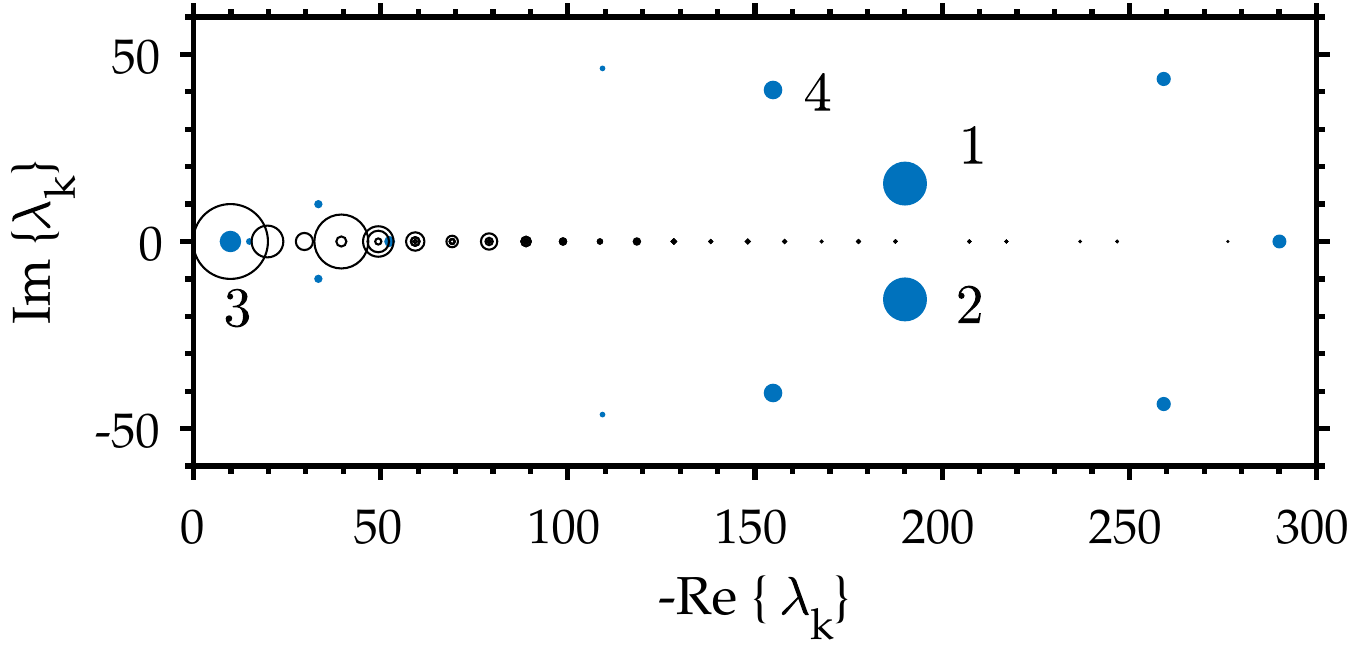}\\
\hspace{1.7cm} a) \hspace{8.5cm} b)
% \includegraphics[width=13cm]{DMD_MODES_Cosine.pdf}\\
% c)\hspace{2cm}
\caption{Results from a DMD analysis of the Burgers' flow $u(t,x)$. Figure (a) shows the accuracy of the DMD. Figure (b) shows the computed DMD eigenvalues. Figure (c) shows the spatial structure of the leading DMD modes.}
\label{DMD_Migu}
\end{figure}

\textcolor{black}{
The DMD solves the regression problem of identifying the eigendecomposition of the best linear propagator $\mathbf{P}$ such that $\mathbf{X}_2=\mathbf{P}\mathbf{X_1}$, with $\mathbf{X}_1,\mathbf{X}_2\in \mathbb{R}^{1023\times 100}$ containing the columns of $\mathbf{X}$ from $1$ to $n_t-1$ and from $2$ to $n_t$ respectively.
The hoped for approximation of the Koopman eigenvalues are the eigenvalues of a reduced propagator $\mathbf{S}=\mathbf{U}^T \mathbf{P}\mathbf{U}$, with $\mathbf{U}$ the matrix of left singular vectors from the singular value decomposition $\mathbf{X}_1=\mathbf{U}\mathbf{\Sigma}\mathbf{V}^T$. Given $\mathbf{S}=\mathbf{W}\mathbf{\Lambda}\mathbf{W}^{-1}$ the eigenvalue decomposition of $\mathbf{S}$, the approximation of the Koopman modes are computed as the columns of the matrix $\mathbf{\Phi}=\mathbf{X}_2 \mathbf{V}\mathbf{\Sigma}^{-1}\mathbf{W}$ as proposed in \cite{Tu2014}. }

\textcolor{black}{
The main results of the DMD analysis are shown in Figure \ref{DMD_Migu}. Figure \ref{DMD_Migu}a) compares the dynamics $u(t,x)$ with the DMD reconstruction using $\mbox{rank}(\mathbf{X})=15$ modes, showing the excellent convergence of the DMD. The associated DMD eigenvalues are shown in \ref{DMD_Migu}(b) (full markers) and compared to the actual Koopman eigenvalues (empty markers). In this figure, the size of the markers is proportional to the mode contribution, computed from formula \eqref{Eq_Importance}. }

\textcolor{black}{
With the only exception of the eigenvalue for $\nu=(1)$, $\lambda_\nu=-\pi^2$, the DMD eigenvalues differ significantly from the actual Koopman ones, with the leading DMD modes (marked as 1 and 2 in the figure) having non-negligible imaginary parts. Consequently, the coefficients and the spatial structures of the DMD differ from those of the Koopman decomposition. 
A detailed analysis of the reasons for such discrepancy is out of the scope of this work and will be presented in a later contribution. In general, an accurate approximation of an operator leads to an accurate approximation of its eigenvalues of multiplicity one, but this does not extend to eigenvalues with higher multiplicity. This example gives a practical illustration of the difficulties one faces in the data-driven identification of the Koopman modes.
}

\section{Concluding Remarks}\label{SecVII}
\textcolor{black}{
	This paper gives an explicit Koopman decomposition (formula \eqref{F3}) for the Burgers' flow. It identifies the coefficients ($\varphi_\nu(u_0)$ in formula \eqref{Phi_nu}) of the decomposition as values of the eigenfunctionals of the Koopman operator (formula \eqref{Ffi}) and the temporal evolution as its eigenvalues. It also identifies the Koopman modes as the coefficients of the decomposition of the Dirac functionals on the Koopman eigenfunctionals (formula \eqref{F4bis}).
	\vskip 0.1cm}
\textcolor{black}{
	This work builds on that of Page and Kerswell \cite{Page2018} and overcomes the intrinsic multiplicity issue they encounter: we construct as many Koopman modes as the geometric multiplicity}. 

\textcolor{black}{The convergence of the Koopman decomposition is here proven. To authors' knowledge, this is the first proof and explicit computation of Koopman modes for the flow given by a PDE.}  % To be polished
\vskip 0.1cm
This paper highlights the need for localization in the function space of the state variable and shows that the convergence of the decomposition is a local property. One must thus restrict the decomposition to its convergence set, which should contain the considered orbits for all $t>t_0$.
In the case of the Burgers equation, the only invariant set, where the decomposition converges, is the vicinity of the sink (the zero solution). This is where the necessary smallness condition shows up, quantified by the definition $\Omega_{\mathcal{B}}$ in \eqref{Omega_0}. For an attempt to localize in a neighborhood of an unstable singular location, an interesting numerical experiment is made by Page and Kerswell in \cite{Page2019}.
\vskip 0.1cm
\textcolor{black}{The need for regularity of the initial condition to have convergence of the Koopman decomposition at $t =0$ is a classic PDE ingredient. It is used in this paper to prove that the Dirac functionals can be decomposed on the eigenfunctionals of the Koopman operator, hence to prove \eqref{phi} for the case $\phi:=\delta_x$.}
\vskip 0.1cm

%%%%%%%%%%%%%%%%%%%%%%%%%%%%   Before %%%%%%%%%%%%%%%%%%%%%%%%%%%%
% Problems I see: 
% 1. starting with 'Here' or with 'So one has to'
% 2. I would remove 'More precisely and link it to one sentence' (this let's us remove 'It'-- there are many of those).
% 3. The sentence with 'restrict to convergence set.Moreover this set should' sounds a bit strange; I would merge them.
% 4. The considered orbits for all $t>t_0$, for some $t_0$ is really strange (if t_0 can be arbitrarily small, then this should just be t>0)

\textcolor{black}{An illustration of the Koopman decomposition is shown. This illustration gives a practical demonstration of the key features of this decomposition, including the multiplicity of the eigenvalues and the difficulties in convergence for $t\rightarrow 0$. These features make the data-driven identification of the actual Koopman modes challenging. To further illustrate this, a classic DMD of the Burgers flow is presented and the DMD eigenvalues compared to the actual Koopman eigenvalues. While the DMD outperforms the Koopman decomposition in terms of convergence and shows no difficulties as $t\rightarrow 0$, only one of the DMD eigenvalues correspond to a Koopman eigenvalue: the leading eigenvalue of multiplicity one.}
\vskip 0.1cm
\textcolor{black}{The explicit Koopman decomposition presented in this work enables a systematic analysis of the feasibility of a data-driven identification of Koopman modes via DMD.}

\section{Appendix: Estimates}\label{SecV}

\subsection {Three Preliminary Properties}\label{lemmi}

Here is the estimate used \textcolor{black}{to derive (\ref{Fdev})}:

{\bf Property 1}: If $v_0$ fulfills $\int_0^1v_0(s)ds=1$ and $\,\Vert \partial_xv_0\Vert_{L^2}<{1\over 4}$  then 
$$\sup_x\vert 1-v_0(x)\vert< {1\over 4}\quad {\rm and}  \quad \Vert u_0 \Vert_{L^2}<{2\over 3}\quad {\rm for} \quad u_0=C(v_0).$$

\noindent {\it Proof:} $v_0$ is given by:
$$v_0(x)=1+\int_0^xs\partial_xv_0(s)ds-\int_x^1(1-s)\partial_xv_0(s)ds
$$
By Cauchy-Schwartz, $\sup_x\vert1-v_0(x)\vert \leq \Vert \partial_xv_0\Vert_{L^2}<{1\over 4}$ so $\vert u_0(x)\vert\leq {8\over 3}\vert \partial_xv_0(x)\vert$.

\noindent
\textcolor{black}{We now give the assumption needed on $u_0 = C(v_0)$ in order that $v_0$ fulfills assumptions of property 1. This assumption on $u_0$ leads to the definition of $\Omega_{\mathcal{B}}$ given by formula (\ref{Omega_0}) that is needed in the proof of formula (\ref{F3}) below (section \ref{proof2}).}

{\bf Property 2:}
Let $u_0 \in L^2$ and $v_0=H(u_0)$ then
$$\sup_x\vert v_0(x)\vert\leq e^{\Vert u_0\Vert_{L^2}} \quad {\rm and}\quad \Vert \partial_xv_0\Vert_{L^2}\leq {1\over 2}e^{ \Vert u_0\Vert _{L^2}}\Vert u_0\Vert_{L^2}$$

\noindent {\it proof:} because $\int_0^xu_0(s)ds\leq \Vert u_0\Vert_{L^2}$ one has $v_0(x)\leq e^{{1\over 2} \Vert u_0\Vert _{L^2}}e^{-{1\over 2} \int_0^xu_0(s)ds}\leq e^{ \Vert u_0\Vert _{L^2}}$. From $u_0v_0=-2\partial_xv_0$ follows  ${ \Vert \partial_xv_0\Vert _{L^2}}\leq {1\over 2}\sup_x\vert v_0(x)\vert  \Vert u_0\Vert _{L^2}\leq  {1\over 2}e^{ \Vert u_0\Vert _{L^2}}\Vert u_0\Vert _{L^2}$.
\\

\noindent \textcolor{black}{Here is the estimate needed to prove uniform and absolute convergence in formula (\ref{F3}) for all $t\ge 0$ for regular initial conditions. It is needed below in section \ref{proof3} to prove formula \eqref{F4}.}

{\bf Property 3:}
Let $u_0\in L^2$, $\partial_xu_0\in L^2$ and $u_0(0)=u_0(1)=0$. If  $v_0=H(u_0)$ then
$$\Vert \partial^2_{xx}v_0\Vert_{L^2}\leq {1\over 2}e^{ \Vert u_0\Vert _{L^2}}\Vert \partial_{x}u_0\Vert_{L^2}(1+{\Vert u_0\Vert_{L^2}\over 2})$$

{\it Proof}: boundary conditions on $u_0$ give $\sup\vert u_0(x)\vert\leq \Vert \partial_{x}u_0\Vert_{L^2}$. 

\noindent Because
$\partial^2_{xx}v_0=-{1\over 2}(v_0\partial_xu_0+u_0\partial_xv_0)$ we get $\Vert \partial^2_{xx}v_0\Vert_{L^2}\leq {1\over 2}(\sup_x\vert v_0(x)\vert \Vert \partial_xu_0\Vert_{L^2} $ $ +\sup_x\vert u_0(x)\vert \Vert \partial_xv_0\Vert_{L^2}) $ and the result follows from the estimates given in property 2.

\subsection {Convergence of formula (\ref{F3})}\label{proof2}
\medskip

\noindent We have, through integration by parts, for all $n\geq 1$: 
$$ l_n(u_0) =  \sqrt{2}\int_0^1H(u_0)(s)\cos{(n\pi s)}ds = 
{\sqrt{2}\over n\pi}\int_0^1\partial_xH(u_0)(s)\sin{(n\pi s)}ds$$
so, because $\Vert u_0 \Vert_{L^2}\leq 1$ for $u_0\in \Omega_{\mathcal{B}}$, we get through Parseval formula and \textcolor{black}{property 2}:
$$\sum_1^\infty n^2\pi^2\vert l_n(u_0)\vert^2=\Vert \partial_xH(u_0)\Vert_{L^2}^2=\Vert \partial_xv_0\Vert_{L^2}^2
\leq {1\over 4}e^{2\Vert u_0\Vert_{L^2}}\Vert u_0\Vert_{L^2}^2\leq {e^2\over 4}\Vert u_0\Vert_{L^2}^2$$
We now prove the absolute convergence of formula (\ref{F3}): because
$$\forall \nu \in \mathcal{A} , \quad\vert a_\nu(x)\vert \leq 2^{\alpha(\nu)+3\over 2}\pi n_0$$
$$I:=\sum_{\nu\in \mathcal{A}}e^{\lambda_{\nu}t}\vert a_\nu(x)\vert \varphi_\nu(u_0)\vert
\leq
2^{3\over 2}\pi\sum_{m=0}^\infty 2^{m\over 2}\sum_{\nu\in \mathcal{A};\alpha(\nu)=m}\prod_{k=0}^me^{-n_k^2\pi^2 t}n_0\prod_{k=0}^m\vert l_{n_k}(u_0)\vert
$$
$$\leq  2^{3\over 2}\pi \,\sum_{n=1}^\infty e^{-n^2\pi^2t}n\vert l_n(u_0)\vert\,\,\sum_{m=0}^\infty 2^{m\over 2}(\sum_{k=1}^\infty \vert l_k(u_0)\vert)^m
$$
We apply the discrete Cauchy-Schwartz inequality and get:

$$I  \leq 2^{3\over 2}\pi \sqrt{\sum_{n=1}^\infty  e^{-2n^2\pi^2t} }\sqrt{\sum_{n=1}^\infty  n^2\vert l_n(u_0)\vert^2 }
\,\,  \sum_{m=0}^\infty 2^{m\over 2}   \beta_2^m  \,\,(\sum_{n=1}^\infty n^2 \vert l_n(u_0)\vert^2 ) ^{m\over 2}
$$
\noindent we use the estimate $\sum n^2\vert l_n(u_0)\vert^2\leq {e^2\over 4\pi^2}\Vert u_0\Vert^2$ proved above and get:
$$I \leq {t^{-{1\over 4}}} e\Vert u_0\Vert_{L^2}\sum_{m=0}^\infty (2^{-{1\over 2}} \pi^{-1}e\beta_2\Vert u_0\Vert_{L^2})^m
$$
\noindent with  $\beta_2=\sqrt{\sum_1^\infty {1/ n^2}}={\pi /\sqrt{6}}$. This series converge because $u_0\in \Omega_{\mathcal{B}}$ so $\Vert u_0\Vert_{L^2}<{ 2\sqrt{3}/ e}$.

We proved convergence of absolute values, hence commutative convergence, as well as uniform convergence in the $x$ variable, \textcolor{black}{for $t >0$.}

\subsection{Uniform convergence in formula (\ref{F4})} \label{proof3}

One only needs to prove $(t,x)$-uniform  convergence  of absolute values for all $t\geq 0, x\in [0,1]$. Let $v_0=H(u_0)$. $v_0$ fulfills Neumann boundary conditions because $u_0 \in  \omega_{\mathcal{B}}$ fulfills Dirichlet boundary conditions.  Moreover $u_0$ has a square integrable weak derivative, so  first and second weak derivatives of $v_0$ are square integrable as shown by property 3. Two integration by parts give:
$$\sum_1^\infty n^4\pi^4\vert l_n(u_0)\vert^2 = \Vert \partial_{xx}^2v_0\Vert_{L^2}^2$$
\textcolor{black}{The estimate goes as follows:}
$$J:=\sum_{\nu\in \mathcal{A}}e^{\lambda_{\nu}t}\vert a_\nu(x)\vert \varphi_\nu(u_0)\vert \leq
2^{3\over 2}\pi\sum_{m=0}^\infty 2^{m\over 2}\sum_{\nu\in \mathcal{A};\alpha(\nu)=m}n_0\prod_{k=0}^m\vert l_{n_k}(u_0\vert)\leq
$$
$$  2^{3\over 2}\pi\sum_{m=0}^\infty 2^{m\over 2}(\sum_{n=1}^\infty n\vert l_n(u_0)\vert)\,\,(\sum_{k=1}^\infty \vert l_k(u_0)\vert)^m
$$
We use the discrete Cauchy-Schwartz inequality to get:
$$J\leq {2\sqrt{2}\over \pi}\beta_2\sqrt{\sum_{n=1}^\infty n^4\vert l_n\vert^2}\sum_{m=0}^\infty (\sqrt{2}\beta_2\sqrt{\sum_{n=1}^\infty n^2\vert l_n\vert^2}\,\,)^m \leq {2^{3\over 2}\beta_2\over \pi}\Vert \partial^2_{xx}v_0\Vert_{L^2}\sum_{m=0}^\infty ({\beta_2e\over\sqrt{2}\pi}\Vert u_0\Vert_{L^2})^m
$$
\noindent this is because $\Vert u_0\Vert_{L^2}\leq 1$ for  $u_0\in \Omega_{\mathcal{B}}$. We use property 3 to get
$$J\leq {3e\beta_2\over \sqrt{2}\pi}\Vert \partial_{x}u_0\Vert_{L^2}\sum_{m=0}^\infty ({\beta_2e\over\sqrt{2}\pi}\Vert u_0\Vert_{L^2})^m
$$
This series converges for $\Vert u_0\Vert_{L^2}<{2\sqrt{3}/ e}$. This condition is fulfilled by any $u_0\in\Omega_{\mathcal{B}}$.

\bibliography{Balabane_Mendez_Najem_2020_arxiv}% Produces the bibliography via BibTeX.

\end{document}